\documentclass[10pt,reqno,final]{amsart}

\usepackage{amsmath,amssymb,amsthm}
\usepackage[top=3.5cm,bottom=3cm,left=3cm,right=3cm]{geometry}
\usepackage{url}

\usepackage[
            pagebackref  = false,
            pdfauthor    = {Giacomo Cherubini and Alberto Perelli},
            pdftitle     = {A spectral universality theorem for Maass L-functions},
            pdfkeywords  = {universality {Maass forms} L-functions},
            pdfcreator   = {Pdflatex},
            pdfpagemode  = UseNone,
            pdfstartview = FitH
           ]
           {hyperref}

\hypersetup{colorlinks=true,pdfborder={0 0 0}}

\numberwithin{equation}{section}

\theoremstyle{theorem}
\newtheorem{theorem}{Theorem}[section]
\newtheorem{proposition}{Proposition}[section]

\newtheorem{lemma}{Lemma}[section]

\theoremstyle{definition}
\newtheorem{remark}{Remark}[section]

\def \ep {\epsilon}

\def \si {\sigma}

\def \d {\mathrm{d}}

\def \bfx {{\boldsymbol{x}}}

\def \CC {\mathbb C}
\def \RR {\mathbb R}
\def \ZZ {\mathbb Z}

\def \A  {{\mathcal A}}
\def \B  {{\mathcal B}}
\def \C  {{\mathcal C}}
\def \D  {{\mathcal D}}
\def \K  {{\mathcal K}}
\def \X  {{\mathcal X}}

\def\fine{\qed}

\newcommand{\GL}{\mathrm{GL}}

\begin{document}

\title{A spectral universality theorem for Maass $L$-functions}
\author{G. Cherubini \lowercase{and} A. Perelli}
\maketitle

\noindent{\bf Abstract.}
We show that for a positive proportion of Laplace eigenvalues $\lambda_j$,
the associated Hecke--Maass $L$-functions $L(s,u_j)$ approximate with
arbitrary precision any target function $f(s)$ on a closed disc with center
in $3/4$ and radius $r<1/4$.
The main ingredients in the proof are the spectral large sieve of
Deshouillers--Iwaniec and Sarnak's equidistribution theorem for Hecke eigenvalues.




\section{Introduction}\label{section1}

A version of Voronin's~\cite{Vor/1975} celebrated universality theorem
for the Riemann zeta function $\zeta(s)$ states that,
given a closed disc $\K$ with radius $0<r<1/4$ centered at $s = 3/4$
and a function $f(s)$ holomorphic and non-vanishing on $\K$,
for any given $\ep>0$ we have
\begin{equation}\label{intro:voronin}
\liminf_{T\to\infty} \frac{1}{T} \, \text{meas}\{ \tau\in [T,2T] : \max_{s\in\K} |\zeta(s+i\tau) - f(s)| < \ep \} >0;
\end{equation}
here ``meas'' denotes the Lebesgue measure.
Thanks to the contribution of many authors, Voronin's theorem has been extended
in various directions, relaxing the conditions on $\K$ and $f$,
or dealing with other $L$-functions in place of $\zeta(s)$
and with various forms of joint universality for vectors of $L$-functions.
We refer to Matsumoto's survey~\cite{Mat/2015} for a fairly complete collection
of such results and several other variants.
Note that universality results for $L$-functions  analogous to \eqref{intro:voronin} do not hold
on the left of the critical line $\si=\Re(s)=1/2$, due to the restrictions
imposed by the functional equation, nor on the right of the
line of absolute convergence $\si=1$, due to boundedness properties;
see the discussion in Perelli--Righetti~\cite{Pe-Ri/rigidity} for further information.
Interestingly enough, there exists a different type of universality that hold also on the right of $\si=1$,
see Andersson--S\"odergren \cite[Theorem 1.10]{An-So/2015}.

We refer to all the above mentioned universality results by the generic name
of universality theorems for families.
In broad terms, a universality theorem for a family of $L$-functions consists
in showing that a positive proportion of the members of the family are arbitrarily close,
uniformly over a given compact set inside the right half of the critical strip,
to any given member of a wide and generic set of target functions.

Besides the family of vertical shifts considered in Voronin's theorem
(and extended to other $L$-functions as observed shortly in the beginning),
there are genuinely GL(1) results for the family of Dirichlet characters,
see Gonek \cite{Gon/1979} and Bagchi \cite{Bag/1981},
or the family of quadratic characters, see Mishou--Nagoshi \cite{MiNa/2006}.
More recently, Kowalski~\cite{Kow/2017} proved a universality theorem for modular $L$-functions in the level aspect,
which can therefore be considered as a genuinely GL(2) statement of this type.
Moreover, he uses the probabilistic setting to give a conceptual explanation
of universality theorems for families of $L$-functions.

In this paper we prove a universality theorem in the eigenvalue aspect for the family
of Hecke--Maass $L$-functions, providing thus a new example of a GL(2)-type result.
Let $\{u_j\}$ be an orthonormal basis of the
space of Hecke-Maass cusp forms for the full modular group,
and let $L(s,u_j)$ be the $L$-function associated with $u_j$. Each $u_j$ satisfies
\[
\Delta u_j = \lambda_j u_j, \quad \lambda_j=\frac{1}{4} + r_j^2 \ \ \text{with} \ \ r_j>0,
\]
and the Weyl law asserts the asymptotic formula
\begin{equation}
\label{1-1}
\sum_{r_j\leq R}1 \sim cR^2
\end{equation}
with $c=1/12$; see e.g.~\cite[Chapter 5]{Iw-Ko/2004}.
Here and in the sequel the $r_j$ are counted with multiplicity,
i.e. each $r_j$ is counted with progressive index as many times
as the (finite) dimension of the eigenspace of $\lambda_j$.
In this way we can fix a one-to-one correspondence between the $r_j$ and the $u_j$,
which is used throughout the paper.
We briefly recall some relevant points of the theory of the Maass forms
and $L$-functions at the beginning of Section~\ref{section2}.

Finally, let $\K$ denote the closed disc $|s-3/4|\leq r$ with $0<r<1/4$,
and let $f(s)$ be a given target function with the following properties:
$f(s)$ is holomorphic and non-vanishing on $\K$, and $f(s)$ is real and positive on $\K\cap\RR$.

\begin{theorem}\label{theorem}
Let $L(s,u_j)$, $\K$ and $f(s)$ be as above. Then, for every $\ep>0$ we have
\[
\liminf_{R\to\infty} \frac{1}{R^2} \# \{r_j\leq R: \max_{s\in\K} |L(s,u_j)-f(s)|<\ep\} >0.
\]
\end{theorem}

\noindent
In other words, in view of~\eqref{1-1} the theorem implies that, given $\ep$, $\K$ and $f(s)$, then as $R\to\infty$ the approximation
\[ 
\max_{s\in\K} |L(s,u_j)-f(s)|<\ep
\]
holds for the forms $u_j$ associated with a (small) positive proportion of $r_j\leq R$.
The reason for choosing the class of target functions $f(s)$ as described in the theorem
can loosely be reconduced to the fact that we need $\log f(s)$ to be well-defined and, since
$\log L(s,u_j)$ has real coefficients, we also need $\log f(s)$ to commute with the complex conjugation.

\begin{remark}\label{rmk-1.1}
It is worth stressing the differences with the situation in the holomorphic case.
It is well known that, for holomorphic modular forms, the Hecke eigenvalues satisfy the Ramanujan conjecture;
moreover, they are distributed according to the Sato--Tate measure.
In the non-holomorphic case, on the other hand, neither of the two corresponding statements is known.
To overcome this fact, we use
the spectral large sieve of Deshouillers--Iwaniec~\cite{De-Iw/1982}, a result of Iwaniec~\cite{Iwa/1992}
that provides the Ramanujan conjecture on average,
and Sarnak's equidistribution theorem for Hecke eigenvalues~\cite{Sar/1987},
which is an analogue of the Sato--Tate conjecture as well as a different average form of the Ramanujan conjecture. \fine
\end{remark}

\begin{remark}\label{rmk-1.2}
The fundamental tools in Voronin's result are the classical Kronecker--Weyl
equidistribution theorem and Pechersky's rearrangement theorem.
Indeed, for vertical shifts of a given $L$-function,
or a given vector of $L$-functions, suitable variations of the Kronecker--Weyl theorem are a typical ingredient.
For different families of $L$-functions, as already noticed in Remark~\ref{rmk-1.1},
different and often sophisticated equidistribution theorems,
peculiar to the family under consideration, are required. 
In fact, in addition to a variety of equidistribution theorems, the techniques in the proof of universality theorems have evolved and enriched since Voronin's first result, for example bringing into play the language and methods of probability theory
and certain density theorems for spaces of holomorphic functions. 
However, one may trace a typical protocol in such proofs. Roughly, in the classical language the protocol leading to a universality result with target functions $f(s)$  on a compact $\K$ for a family of $L$-functions $L(s,u)$, with $u$ in a set $U$ of size $|U|$, may be outlined as follows.

\vskip-.12cm
(i) Use suitable averages over $U$ (often depending on $L^2$-bounds) to replace $L(s,u)$, or $\log L(s,u)$, by finite approximations, uniformly for $s\in\K$ and for all $u$ in a subset $A\subset U$ of (normalized) measure 1, as $|U| \to \infty$; $\log L(s,u)$ is employed in the case of Euler products, since usually the equidistribution theorems involve prime numbers.

\vskip-.12cm
(ii) Use suitable density theorems in spaces of holomorphic functions to approximate,
uniformly over $\K$, the target functions $f(s)$, or their logarithm, by a certain class of random Dirichlet polynomials.

\vskip-.12cm
(iii) Use suitable equidistribution theorems to show that the approximations in (i) become, as $|U| \to \infty$, arbitrarily close to the random Dirichlet polynomials in (ii), uniformly over $\K$ and for all $u$ in a subset of $A$ of (small) positive measure.

\vskip-.12cm
\noindent
The above triple approximation procedure leads the required universality result,
and can be traced in the probabilistic setting as well; see in particular Kowalski~\cite{Kow/2017}.
We have chosen to present our result in the classical language;
our approach might therefore be conceptually less satisfactory, but probably is more transparent. \fine
\end{remark}

As an example, Voronin's theorem deals with a single $L$-function, namely $\zeta(s)$,
and the parameters $u\in U$ are the vertical shifts $\tau\in[T,2T]$.
Then step (i) follows from standard mean square approximations of $\zeta(s)$ by finite Euler products,
step (ii) follows from Pechersky's theorem, and step (iii) is a consequence of the Kronecker--Weyl theorem.
In our case, step (i) follows from the spectral large sieve of Deshouillers--Iwaniec~\cite{De-Iw/1982},
whereas step (iii) is performed by applying Sarnak's equidistribution theorem~\cite{Sar/1987};
step (ii) depends in a rather mild way on the $L$-functions into play, and is essentially
independent from the rest of the paper.

We outline the proof of the theorem in Section~\ref{section2}; details are then given in Section~\ref{section3}. Finally,
a standard argument allows to relax the requirements on $f(s)$ slightly, namely replacing holomorphy and non-vanishing on $\K$ by holomorphy and non-vanishing on the interior of $\K$ and continuity on the boundary of $\K$. Moreover, we expect that with more work one could choose more general compact sets $\K$ in the semistrip $1/2<\si<1$.

{\bf Acknowledgements.}
G.C. was supported by the INdAM grant {\sl Ing.G.Schirillo} 2017--2018;
A.P. is member of the INdAM group GNAMPA. We thank A. S\"odergren for bringing \cite{An-So/2015} to our attention.

\section{Outline of the proof}\label{section2}

As anticipated in the introduction, we first review some properties of the $u_j$ and $L(s,u_j)$,
then we gather the key ingredients to prove steps (i)--(iii) of the protocol described in Remark~\ref{rmk-1.2},
and finally we show how Theorem~\ref{theorem} follows from such ingredients.

\subsection{Hecke--Maass forms and $L$-functions}

Let us start with reviewing the basic notions concerning
Hecke--Maass cusp forms required for our needs.
A reference for this material is Chapter 3 of Goldfeld's book~\cite{Gol/2006}.
Let $\{u_j\}$, $j\geq 1$, be an orthonormal basis of Hecke--Maass cusp forms of weight 0
for the full modular group PSL$(2,\ZZ)$.
These are simultaneous eigenfunctions of the hyperbolic Laplacian $\Delta$
and of the Hecke operators $T_p$, therefore we can write
\[
\begin{cases}
T_pu_j =\lambda_j(p) u_j \\
\Delta u_j = \lambda_j u_j = (\frac{1}{4} + r_j^2) u_j.
\end{cases}
\]
Recall that $\lambda_j(p), \lambda_j\in\RR$, and since $\lambda_j>1/4$ we choose $r_j> 0$.
Set $n(p)$ to be the operator norm of $T_p$, which equals
\begin{equation}
\label{2-1}
n(p) = p^{1/2} + p^{-1/2}.
\end{equation}
Clearly $|\lambda_j(p)|\leq n(p)$, and the Ramanujan conjecture predicts that for $\lambda_j\neq0$ we should have the much stronger bound $|\lambda_j(p)|\leq 2$ (note that $n(p)>2$ for all $p$). The best known result in this direction is due to Kim--Sarnak, see~\cite{Kim/2003}, who proved that $\lambda_j(p) \ll p^{7/64}$. We also recall that every $u_j$ has a Fourier--Bessel expansion. Indeed, writing $z=x+iy$ and $e(x)=e^{2\pi ix}$, we have
\[
u_j(z) = y^{1/2} \sum_{n\neq0} \rho_j(n) K_{ir_j}(2\pi|n|y)e(nx),
\]
where $K_\nu$ is the $K$-Bessel function and the coefficients $\rho_j(n)$ satisfy $\rho_j(-n) = \overline{\rho_j(n)}$ and $\rho_j(1)\neq0$ for all $j$.
The Hecke--Maass $L$-function associated with $u_j$ is defined for $\si>1$ as
\begin{equation}
\label{2-2}
L(s,u_j) = \prod_p \left(1-\frac{\lambda_j(p)}{p^s} + \frac{1}{p^{2s}}\right)^{-1} = \sum_{n=1}^\infty \frac{\lambda_j(n)}{n^s}.
\end{equation}
It extends to $\CC$ as an entire function of finite order with polynomial growth on vertical strips and satisfies the  functional equation
\[
\Lambda(s,u_j) = (-1)^{\ep_j} \Lambda(1-s,u_j),
\]
where $\ep_j$ equals 0 or 1 depending on whether $u_j$ is even or odd (i.e. $u_j(-x+iy) = \pm u_j(x+iy)$) and
\begin{equation}
\label{2-3}
\Lambda(s,u_j) = \pi^{-s} \Gamma\left(\frac{s+\ep_j+ir_j}{2}\right) \Gamma\left(\frac{s+\ep_j-ir_j}{2}\right) L(s,u_j) =  \gamma_j(s) L(s,u_j),
\end{equation}
say. For every $r_j$ we have $\lambda_j(1)=1$, and for $n\geq1$ the coefficients
$\lambda_j(n)$ and $\rho_j(n)$ are related by
\begin{equation}
\label{2-4}
\rho_j(n) = \rho_j(1) \lambda_j(n).
\end{equation}
Moreover, for $\si>1$ we may rewrite the Euler product in~\eqref{2-2} as
\begin{equation}
\label{2-5}
L(s,u_j) = \prod_p \Big(1-\frac{\alpha_j(p)}{p^s}\Big)^{-1} \Big(1-\frac{\beta_j(p)}{p^s}\Big)^{-1},
\end{equation}
and it is known that $\alpha_j(p),\beta_j(p)\in\RR$, with $\alpha_j(p)\beta_j(p)= 1$
and $\alpha_j(p) + \beta_j(p) = \lambda_j(p)$, which in turn implies by~\eqref{2-1} the inequality
\begin{equation}
\label{2-6}
\max\big(|\alpha_j(p)|, |\beta_j(p)|\big) \leq p^{1/2}
\end{equation}
for every $j$ and $p$.

\subsection{Average bounds.}

The approximations alluded to in step (i) of Remark~\ref{rmk-1.2}
are based on the well-known spectral large sieve by Deshouillers--Iwaniec~\cite{De-Iw/1982},
coupled with other results by Iwaniec~\cite{Iwa/1992} and Luo~\cite{Luo/1999}.

\begin{lemma}[{\cite[Theorem 2]{De-Iw/1982}}]\label{lemma-2.1}
Let $(a_n)$ be any sequence of complex numbers and $\ep>0$. Then
\begin{equation}
\label{2-7}
\sum_{r_j\leq R} \frac{1}{\cosh \pi r_j} \big| \sum_{n\leq N} a_n \rho_j(n) \big|^2 \ll (R^2 + N^{1+\ep}) \sum_{n\leq N} |a_n|^2.
\end{equation}
\end{lemma}

In view of~\eqref{2-4}, the transition from the coefficients $\rho_j(n)$ in~\eqref{2-7}
to the coefficients $\lambda_j(n)$ of $L(s,u_j)$ requires a sharp control on average on the numbers
\begin{equation}
\label{2-8}
\alpha_j = \frac{|\rho_j(1)|^2}{\cosh \pi r_j}.
\end{equation}
It follows from results of~\cite{Iwa/1990} and~\cite{hoffstein-lockhart} that $r_j^{-\ep}\ll \alpha_j \ll r_j^\ep$, and the required average bound is achieved by a weak form of a special case of Theorem 1 in~\cite{Luo/1999}, which we report as follows.
\begin{lemma}[{\cite[Theorem 1]{Luo/1999}}]\label{lemma-2.2}
Let $\alpha_j$ be defined as in~\eqref{2-8}. Then
\[
\sum_{r_j \leq R} \alpha_j^{-1} \ll R^2.
\]
\end{lemma}
We also need the following uniform bound for the $\lambda_j(n)^2$, proved in~\cite{Iwa/1992},
which can be viewed as a Ramanujan conjecture on average for $\lambda_j(n)$.
\begin{lemma}[{\cite[Lemma 1]{Iwa/1992}}]\label{lemma-2.3}
Let $\ep>0$. Then
\[
\sum_{n\leq N} \lambda_j(n)^2 \ll r_j^\ep N.
\]
\end{lemma}
\noindent%
For a given $X>2$ we define the inverse partial product of $L(s,u_j)$ by
\begin{equation}
\label{2-9}
L_X(s,u_j) := \prod_{p\leq X}  \left(1-\frac{\lambda_j(p)}{p^s} + \frac{1}{p^{2s}}\right).
\end{equation}
Note, thanks to~\eqref{2-5} and~\eqref{2-6}, that $L_X(s,u_j)$ is holomorphic and non-vanishing for $\si>1/2$. Given a closed disc $\K$ as in Section~\ref{section1}, i.e. with radius $0<r<1/4$ and center at $s = 3/4$, we define
\begin{equation}
\label{2-10}
\theta := \min_{s\in \K} \Re(s) - \frac{1}{2}.
\end{equation}
In other words, $\theta$ is the distance of $\K$ from the critical line.
Clearly $\theta<1/4$, and we stress that it is crucial for the proof of Theorem~\ref{theorem} having $\theta>0$.
The average bounds that we seek are summarized by the following proposition,
based on Lemmas~\ref{lemma-2.1}--\ref{lemma-2.3}.

\begin{proposition}\label{prop-2.1}
Let $2<Y<X$ be sufficiently large and $R>0$ be much larger than $X$. Then
\begin{equation}
\label{2-11}
\sum_{r_j\leq R} \max_{s\in\K} |L(s,u_j) L_X(s,u_j) -1| \ll R^2 X^{-\theta/2}
\end{equation}
and
\begin{equation}
\label{2-12}
\sum_{r_j\leq R} \max_{s\in\K} |\log L_Y(s,u_j) - \log L_X(s,u_j)| \ll R^2 Y^{-\theta/2}.
\end{equation}
\end{proposition}
By ``much larger'' in Proposition~\ref{prop-2.1}
we mean, for example, that $R \geq e^{e^X}$; actually, in the end we let $R\to \infty$, while $X$ and $Y$ will be large depending on $\ep$. Therefore, to the purpose of proving Theorem~\ref{theorem}, this assumption is not restrictive.
The proof of Proposition~\ref{prop-2.1} is given in Section~\ref{section-3.1}.

\subsection{Density in H-spaces}

Let $H_\RR(\D)$ be the space of holomorphic functions $g(s)$ on the open disc $\D=\{|s-3/4|<r\}$ which are continuous on $\K=\overline{\D}$  and real on $\D\cap \RR$, with the norm
\[
\| g \| = \max_{s\in\K} |g(s)|.
\]
Note that if $f(s)$ is a target function as in Theorem~\ref{theorem}, then
\begin{equation}
\label{2-13}
g(s) = \log f(s) \ \text{belongs to} \ H_\RR(\D). 
\end{equation}
Note also that the requirement that $g(s)$ is real on $\D\cap \RR$ comes from the fact that the Dirichlet series $\log L(s,u_j)$ have real coefficients.
As outlined in (ii) in Section~\ref{section1}, our basic tool in this step is a suitable version of the density results for spaces of holomorphic functions employed, since the early researches of Good~\cite{Goo/1981}, Gonek~\cite{Gon/1979} and, soon after, independently by Bagchi~\cite{Bag/1981}, by various authors in the context of universality theorems. Here we refer to a lemma of Kowalski~\cite{Kow/2017}, which we rewrite as follows.

\begin{lemma}[{\cite[Lemma 6]{Kow/2017}}]\label{lemma-2.4}
For any given $Z> 2$ the set of all sums of type
\[
\sum_{Z<p\leq Y} \frac{\omega_p}{p^s},
\]
with $Y>Z$ and $\omega_p\in[-2,2]$, is dense in $H_\RR(\D)$.
\end{lemma}
Lemma~\ref{lemma-2.4} is an immediate consequence of Lemma 6 in~\cite{Kow/2017},
where the numbers $\omega_p$ are viewed as $\mathrm{Tr}(x_p)$ with $x_p\in\mathrm{SU}_2(\CC)$.
From Lemma \ref{lemma-2.4} we deduce the required density result in the following form.
\begin{proposition}\label{prop-2.2}
Let $g\in H_\RR(\D)$ and $\ep>0$. Then there exist arbitrarily large $Y>Y_\ep>2$ and $\omega_p\in[-2,2]$ for $p\leq Y$ such that
\[
\max_{s\in\K} \Big| g(s) + \sum_{p\leq Y} \log\Big(1-\frac{\omega_p}{p^s} + \frac{1}{p^{2s}}\Big)\Big| <\ep.
\]
\end{proposition}
The simple proof of Proposition~\ref{prop-2.2} is given in Section~\ref{section-3.2}.

\subsection{Equidistribution}\label{subsection-equidistribution}

The last step in the protocol described in Remark~\ref{rmk-1.2}
consists, in our case, in the application of Sarnak's equidistribution theorem in the $j$-aspect
for the coefficients $\lambda_j(p)$ of the Hecke--Maass $L$-functions.

For a given $r_j$ we consider the infinite vector $\bfx_j=(\lambda_j(p))$, $p$ prime. In view of the bound $|\lambda_j(p)|\leq n(p)$, see after~\eqref{2-1}, the vectors $\bfx_j$ belong to the product space
\[
\X = \prod_p \X_p,\quad \X_p=[-n(p),n(p)].
\]

\begin{lemma}[{\cite[Theorem 1.2]{Sar/1987}}]\label{lemma-2.5}
As $j\to\infty$ the vectors $\bfx_j$ become equidistributed in $\X$
with respect to the product measure $\mu = \prod_p\mu_p$, where
\[
\d\mu_p(x)=
\begin{cases}
\displaystyle\frac{1}{2\pi} \frac{(1+p)\sqrt{4-x^2}}{n(p)^2-x^2}\d x & \text{{\it if }} \, |x|\leq 2,\\
0 & \text{{\it otherwise}}.
\end{cases}
\]
\end{lemma}

Note that, despite the endpoints of the intervals $\X_p$ are of size $p^{1/2}$,
Sarnak's theorem implies that, when $j\to\infty$, the eigenvalues $\lambda_j(p)$ tend to concentrate on the smaller interval $[-2,2]$, for all~$p$. Of course, this supports the Ramanujan conjecture. Note also that $\mu_p$ are probability measures, meaning that
\[
\int_{\X_p} \d\mu_p(x) = 1.
\]
Moreover, what really matters for us is not the exact form of the measures $\d\mu_p$ in Lemma~\ref{lemma-2.5}, but just the fact that if $[a,b]\subset [-2,2]$, then $\mu_p([a,b]) >0$ for every $p$. We further recall that Lemma~\ref{lemma-2.5} implies that for any $\phi\in C_c(\X)$ we have
\begin{equation}
\label{2-14}
\lim_{R\to\infty} \frac{1}{cR^2} \sum_{r_j\leq R} \phi(\bfx_j) =  \int_\X \phi(\bfx) \d \mu(\bfx),
\end{equation}
where $c>0$ is the constant in~\eqref{1-1} and $\bfx=(x_p)$.

The equidistribution theorem~\eqref{2-14} is crucial in the proof of the following proposition. Let $Y>2$, $0<\delta<10^{-3}$, $\omega_p\in[-2,2]$ for $p\leq Y$ be given and let
\[
\Omega_p= \Omega_p(\delta) =[\omega_p-\delta, \omega_p+\delta]\subset \X_p,
\]
\[
\B_R = \B_R(\delta,Y) =  \{r_j\leq R : \lambda_j(p)\in\Omega_p \ \text{for every} \ p\leq Y\}.
\]
Clearly, $\Omega_p$ and $\B_R$ may depend also on the $\omega_p$, but such a dependence is less important here since in the end the $\omega_p$ will depend only on the target function $f(s)$.

\begin{proposition}\label{prop-2.3}
Let $R$ be sufficiently large, $c>0$ be as in~\eqref{1-1}, $0<\delta<10^{-3}$  and $X>Y>2$. Then
\begin{equation}
\label{2-15}
|\B_R| \geq \frac{c}{4}  \Big(\prod_{p\leq Y} \mu_p(\Omega_p)\Big) R^2,
\end{equation}
\begin{equation}
\label{2-16}
\max_{s\in\K} \Big| \log L_Y(s,u_j) - \sum_{p\leq Y} \log\Big(1-\frac{\omega_p}{p^s} + \frac{1}{p^{2s}} \Big)\Big| \ll \delta \sqrt{Y}
\end{equation}
for every $r_j\in\B_R$, and
\begin{equation}
\label{2-17}
\begin{split}
\sum_{r_j\in\B_R} &\max_{s\in\K} \Big| \log L_X(s,u_j) - \log L_Y(s,u_j) \Big|  \\
&\leq 4 \Big(\prod_{p\leq Y} \mu_p(\Omega_p)\Big)  \sum_{r_j\leq R} \max_{s\in\K} \Big| \log L_X(s,u_j) - \log L_Y(s,u_j) \Big|.
\end{split}
\end{equation}
\end{proposition}

The proof of Proposition~\ref{prop-2.3} is given in Section~\ref{section-3.3}.

\subsection{Proof of Theorem~\ref{theorem}}

Finally, we show that the theorem follows easily from Propositions~\ref{prop-2.1},
\ref{prop-2.2}, and~\ref{prop-2.3}.

Let $\ep>0$ be arbitrarily small.  From~\eqref{2-11} in Proposition~\ref{prop-2.1} we deduce that
\[
\max_{s\in\K}|L(s,u_j) L_X(s,u_j) -1| <\ep/2
\]
for every $r_j\leq R$ in a subset $\A_R$ of cardinality
\begin{equation}
\label{2-18}
|\A_R| \geq R^2 \big(1-O(X^{-\theta/2}/\ep)\big).
\end{equation}
Since $L_X(s,u_j)$ is holomorphic and non-vanishing for $\si>1/2$, in particular we have that $L(s,u_j)\neq 0$ for all $s\in\K$ and every $r_j\in\A_R$. Hence for such $r_j$ the logarithm of $L(s,u_j)$ (and of $L_X(s,u_j)$) is well defined and holomorphic on $\K$. Moreover, we have that
\begin{equation}
\label{2-19}
\max_{s\in\K}|\log L(s,u_j) + \log L_X(s,u_j)| < \ep
\end{equation}
for every $r_j\in\A_R$.

Let now $f(s)$ be a target function as in the theorem. Thanks to \eqref{2-13} we can apply Proposition~\ref{prop-2.2} to $g(s)=\log f(s)$, thus getting arbitrarily large $Y>Y_\ep>2$ and real numbers $\omega_p\in[-2,2]$ for $p\leq Y$ such that
\begin{equation}
\label{2-20}
\max_{s\in\K} \Big| \log f(s) + \sum_{p\leq Y} \log\Big(1-\frac{\omega_p}{p^s} + \frac{1}{p^{2s}}\Big)\Big| <\ep.
\end{equation}

Next, arguing as for \eqref{2-19}, from \eqref{2-15} and \eqref{2-17} of Proposition~\ref{prop-2.3} and \eqref{2-12} of Proposition \ref{prop-2.1} we see that
\begin{equation}
\label{2-21}
 \max_{s\in\K} \Big| \log L_X(s,u_j) - \log L_Y(s,u_j) \Big| <\ep
\end{equation}
for all $r_j$ in a subset $\C_R$ of $\B_R$ of cardinality 
\begin{equation}
\label{2-22}
|\C_R|\geq \frac{1}{2} |\B_R|,
\end{equation}
provided $Y=\overline{Y}_\ep>2$ is sufficiently large and \eqref{2-20} holds. Hence, now we choose $\delta=\delta_\ep>0$ in Proposition~\ref{prop-2.3} so small that, with this choice of $Y$, \eqref{2-16} becomes
\begin{equation}
\label{2-23}
\max_{s\in\K} \Big| \log L_Y(s,u_j) - \sum_{p\leq Y} \log\Big(1-\frac{\omega_p}{p^s} + \frac{1}{p^{2s}} \Big)\Big| <\ep
\end{equation}
for every $r_j\in\C_R$. Moreover, we further choose $X=X_\ep$ so large that, thanks to \eqref{2-15}, \eqref{2-18} and \eqref{2-22}, the subset
\begin{equation}
\label{2-24}
\D_R = \A_R \cap \C_R \quad \text{has cardinality} \quad |\D_R| \gg_\ep R^2.
\end{equation}

Finally, thanks to \eqref{2-24} we have that \eqref{2-19}, \eqref{2-20}, \eqref{2-21} and \eqref{2-23} hold simultaneously for every $r_j\in\D_R$, thus for such $r_j$ we deduce
\begin{equation}
\label{2-25}
\max_{s\in\K} |\log L(s,u_j) - \log f(s)| < 4\ep.
\end{equation}
But from \eqref{2-25} we immediately obtain also that
\begin{equation}
\label{2-26}
\max_{s\in\K} |L(s,u_j) - f(s)| < 5\ep
\end{equation}
for all $r_j\in\D_R$. The theorem follows now from \eqref{2-24} and \eqref{2-26} since $\ep>0$ is arbitrary. \fine

\section{Proof of the propositions}\label{section3}

In this section we always assume that $R$ is much larger than $X$ and $Y$.

\subsection{Proof of Proposition~\ref{prop-2.1}}\label{section-3.1}

We start with a Dirichlet polynomial approximation to $L(s,u_j)$ for $s\in\K$ and $r_j\leq R$. Let $N>1$ be sufficiently large, $\ep>0$, $\varphi:[0,\infty)\to[0,1]$ be a smooth function compactly supported in $[0,2]$ with $\varphi(x)=1$ for $x\in[0,1]$, and let $\widetilde{\varphi}(w)$ be its Mellin transform. Then by a standard technique, see e.g. Chapter 5 of Iwaniec--Kowalski~\cite{Iw-Ko/2004}, for $1/2\leq \si\leq 1$ we obtain that
\begin{equation}
\label{3-1}
L(s,u_j) = \sum_{n\leq 2N} \frac{\lambda_j(n)}{n^s} \varphi\left(\frac{n}{N}\right) + (-1)^{\ep_j}
\sum_{n= 1}^\infty \frac{\lambda_j(n)}{n^{1-s}} W(nN,s),
\end{equation}
where
\begin{equation}
\label{3-2}
W(\xi,s) = - \frac{1}{2\pi i} \int_{(1+\ep)} \frac{\gamma_j(1-s+w)}{\gamma_j(s-w)} \widetilde{\varphi}(-w) \xi^{-w} \d w
\end{equation}
and $\gamma_j(s)$ is as in \eqref{2-3}. The convergence of the integral in \eqref{3-2} and of the second sum in \eqref{3-1} is granted by Stirling's formula and the decay property $\widetilde{\varphi}(w)\ll|w|^{-A}$, for arbitrarily large $A$.
Moreover, using Lemma~\ref{lemma-2.3} we see that the second sum in \eqref{3-1} is $O(r_j^{2-2\theta + 5\ep/2}N^{-1-\ep})$, uniformly for $s\in\K$, $\theta$~being defined by \eqref{2-10}. Hence, choosing 
\begin{equation}
\label{3-3}
N=R^{2-\theta}, 
\end{equation}
we get the desired approximation, i.e. uniformly for $r_j\leq R$ and $s\in\K$ we have
\begin{equation}
\label{3-4}
L(s,u_j) = \sum_{n\leq 2N} \frac{\lambda_j(n)}{n^s} \varphi\left(\frac{n}{N}\right) + O(R^{-\theta+\ep}).
\end{equation}

\noindent%
Next, recalling \eqref{2-9} we observe that
\begin{equation}
\label{3-5}
L_X(s,u_j) = \sum_{\substack{m\leq P(X)\\ p|m\Rightarrow p\leq X}} \sum_{kl^2=m} \frac{\lambda_j(k)\mu(k)|\mu(kl)|}{m^s}
\end{equation}
with 
\begin{equation}
\label{3-6}
P(X) = \prod_{p\leq X} p^2 \ll e^{(2+\ep)X},
\end{equation}
hence bounding trivially and using Lemma~\ref{lemma-2.3} we obtain, for every $s\in\K$, the estimate
\begin{equation}
\label{3-7}
L_X(s,u_j) \ll r_j^{\ep} e^{(1+\ep)X}.
\end{equation}
Therefore, observing from \eqref{3-3} that $N>X$, from \eqref{2-9},\eqref{3-4},\eqref{3-5} and \eqref{3-7} we deduce that
\begin{equation}
\label{3-8}
\begin{split}
L(s,u_j) L_X(s,u_j) &= 1 +   \sum_{X<u\leq 2NP(X)} \frac{1}{u^s} \sum_{\substack{nkl^2=u \\ p|kl^2 \Rightarrow p\leq X \\ n\leq 2N}}  \lambda_j(n) \lambda_j(k)\mu(k)|\mu(kl)|\varphi\left(\frac{n}{N}\right) \\
&\hskip1.5cm + O(R^{-\theta+2\ep} e^{(1+\ep)X}) \\
&= 1 + \Sigma_{X,N}(s,u_j) + O(R^{-\theta+2\ep} e^{(1+\ep)X}),
\end{split}
\end{equation}
say, uniformly for $r_j\leq R$ and $s\in\K$. We wish now to apply Lemma~\ref{lemma-2.1}. In order to do so, we first rearrange the sum $\Sigma_{X,N}(s,u_j)$. Since the coefficients $\lambda_j(n)$ satisfy the multiplication rule of the Hecke operators, namely
\[
\lambda_j(n)\lambda_j(k) = \sum_{d|(n,k)} \lambda_j\left(\frac{nk}{d^2}\right),
\]
we can write
\[
\Sigma_{X,N}(s,u_j) = \sum_{X<u\leq 2NP(X)} \frac{1}{u^s} \sum_{\substack{nkd^2l^2=u \\ p|kdl^2 \Rightarrow p\leq X \\ d,n\leq 2N}}  \lambda_j(nk) \mu(kd)|\mu(kdl)|\varphi\left(\frac{nd}{N}\right).
\]
Hence rearranging we get the estimate
\begin{equation}
\label{3-9}
\Sigma_{X,N}(s,u_j)  \ll  \sum_{d\leq 2N} \frac{1}{d^{2\si}} \sum_{l\leq \frac{ \sqrt{2NP(X)}}{d}} \frac{1}{l^{2\si}} \Big| \sum_{\frac{X}{d^2l^2} < \nu \leq \frac{2N P(X)}{d^2l^2}} \frac{\lambda_j(\nu)}{\nu^s} \beta_{d,l}(\nu) \Big|,
\end{equation}
where clearly
\begin{equation}
\label{3-10}
\beta_{d,l}(\nu) \ll \nu^\ep \qquad \text{uniformly in $d$ and $l$.}
\end{equation}
Now we multiply and divide \eqref{3-9} by $\alpha_j^{1/2}$ (see \eqref{2-8}), take the maximum over $s\in\K$ and sum for $r_j\leq R$. Then by Cauchy's formula, after switching summation and integration we obtain
\[
\sum_{r_j\leq R} \max_{s\in\K} |\Sigma_{X,N}(s,u_j)| \ll \iint_{\widetilde{\K}} \Big(\sum_{r_j\leq R} |\Sigma_{X,N}(s,u_j)|\Big) \d\si \d t,
\]
and by the Cauchy--Schwarz inequality we can further bound the above by
\[
\ll \sum_{d\geq1}\sum_{l\geq1} \frac{1}{(dl)^{1+2\widetilde{\theta}}} \Big(\sum_{r_j\leq R} \alpha_j^{-1}\Big)^{1/2} 
 \Big(\iint_{\widetilde{\K}} \Big( \sum_{r_j\leq R} \frac{1}{\cosh \pi r_j} \Big| \sum_{\frac{X}{d^2l^2} < \nu \leq \frac{2N P(X)}{d^2l^2}} \frac{\rho_j(\nu)}{\nu^s} \beta_{d,l}(\nu) \Big|^2 \Big) \d \si \d t\Big)^{1/2},
\]
where $\widetilde{\K}$ is a suitable open disc containing $\K$,
and $\widetilde{\theta}>0$ is the distance of $\widetilde{\K}$ from the critical line.
Hence by Lemma~\ref{lemma-2.1} and Lemma~\ref{lemma-2.2},
together with \eqref{3-3},\eqref{3-6}, and \eqref{3-10}, we have
\[
\begin{split}
\sum_{r_j\leq R} \max_{s\in\K} |\Sigma_{X,N}(s,u_j)| &\ll R \sum_{d\geq1}\sum_{l\geq1} \frac{1}{(dl)^{1+2\widetilde{\theta}}}  \left(R^2+\left(\frac{NP(X)}{d^2l^2}\right)^{1+\ep}\right)^{1/2} \Bigg(\sum_{v> \frac{X}{d^2l^2}} \frac{1}{v^{1+2\widetilde{\theta}-2\ep}} \Bigg)^{1/2} \\
&\ll \big(R^2+RN^{1/2+\ep/2}e^{(1+\ep/2)X}\big) X^{-\widetilde{\theta}+\ep} \ll R^2 X^{-\widetilde{\theta}+\ep}
\end{split}
\]
since $R$ is much larger than $X$.
The first assertion of Proposition~\ref{prop-2.1} follows now in view of \eqref{1-1} and \eqref{3-8}. The second assertion follows in a very similar way. \fine

\subsection{Proof of Proposition~\ref{prop-2.2}}\label{section-3.2}

Given $\ep>0$ and recalling \eqref{2-10}, we choose an arbitrarily large $Z> Z_\ep\geq2$ such that
\[
2\sum_{p>Z}\sum_{m=2}^\infty \frac{4^m}{p^{m(1/2+\theta)}} <\ep/2
\]
and for $s\in\K$ we write
\[
g_Z(s) = g(s) + \sum_{p\leq Z} \log\Big(1+\frac{1}{p^{2s}}\Big).
\]
By Lemma~\ref{lemma-2.4} there exists $Y=Y_\ep>Z$ (hence $Y$ is also arbitrarily large) and $\omega_p\in[-2,2]$ for $Z<p\leq Y$ such that
\[
\max_{s\in\K} \Big| g_Z(s) - \sum_{Z<p\leq Y} \frac{\omega_p}{p^s}\Big| < \ep/2.
\]
Therefore, letting $\omega_p=0$ for $p\leq Z$, we have
\[
\begin{split}
\max_{s\in\K} \Big| g(s) + \sum_{p\leq Y} \log\Big(1-\frac{\omega_p}{p^s} + \frac{1}{p^{2s}}\big) \Big|
&=
\max_{s\in\K}  \Big| g_Z(s) + \sum_{Z<p\leq Y} \log\Big(1-\frac{\omega_p}{p^s} + \frac{1}{p^{2s}}\big) \Big|
\\
& \leq
\max_{s\in\K} \Big| g_Z(s) - \sum_{Z<p\leq Y} \frac{\omega_p}{p^s}\Big| + 2\sum_{Z<p\leq Y} \sum_{m=2}^\infty  \frac{4^m}{p^{m(1/2+\theta)}} < \ep,
\end{split}
\]
and the result follows. \fine

\subsection{Proof of Proposition~\ref{prop-2.3}}\label{section-3.3}

Recall, for $p\leq Y$, the definition of $\Omega_p$ in Section~\ref{subsection-equidistribution}, and write
\[
\widetilde{\Omega}_p =[\omega_p-\delta-\delta^2, \omega_p-\delta] \cup [\omega_p+\delta, \omega_p+\delta+\delta^2].
\]
Choosing a function $\phi\in C_c(\X)$, depending only on the variables $x_p$ with $p\leq Y$ and satisfying
\[
\phi(\bfx) = 1 \ \text{if} \  x_p\in\Omega_p, \ \ 0\leq \phi(\bfx) \leq 1 \ \text{if} \  x_p\in\ \widetilde{\Omega}_p, \ \ \phi(\bfx) = 0 \ \text{if} \ x_p\in\X_p\setminus (\Omega_p\cup\widetilde{\Omega}_p),
\]
from \eqref{2-14} we obtain that for $R$ sufficiently large and $\delta$ sufficiently small the cardinality of $\B_R$ satisfies 
\[
|\B_R| \geq \frac{c}{4}  \Big(\prod_{p\leq Y} \mu_p(\Omega_p)\Big) R^2
\]
with $c$ as in \eqref{1-1}, and the first assertion follows. Moreover, for every $r_j\in\B_R$ we have
\[
\max_{s\in\K} \Big| \log L_Y(s,u_j) - \sum_{p\leq Y} \log\Big(1-\frac{\omega_p}{p^s} + \frac{1}{p^{2s}} \Big)\Big| \ll  \max_{s\in\K} \Big|\sum_{p\leq Y} \frac{\delta}{p^s} \Big| \ll \delta \sqrt{Y},
\]
thus proving the second assertion.

The third assertion of Proposition~\ref{prop-2.3} follows by a double application of the equidistribution theorem \eqref{2-14}. First we choose the function $\psi(\bfx)$ depending only on the variables $x_p$ with $p\leq X$ defined by
\[
\psi(\bfx) = 
\begin{cases}
\phi(\bfx) & \text{for $p\leq Y$,} \\
\displaystyle\max_{s\in\K}\Big| \sum_{Y<p\leq X} \log\Big(1-\frac{x_p}{p^s} + \frac{1}{p^{2s}}\Big) \Big| \ \text{if} \ x_p\in\X_p  & \text{for $Y<p\leq X$}.
\end{cases}
\]
Then by \eqref{2-14} we have that for $\delta$ sufficiently small
\begin{equation}
\label{3-11}
\begin{split}
\lim_{R\to\infty} & \frac{1}{cR^2} \sum_{r_j\in\B_R} \max_{s\in\K} \Big| \sum_{Y<p\leq X} \log\Big(1-\frac{\lambda_j(p)}{p^s} + \frac{1}{p^{2s}}\Big) \Big| \\
&= \lim_{R\to\infty} \frac{1}{cR^2} \sum_{r_j\leq R} \psi(\bfx_j)\\
&\leq 2  \lim_{R\to\infty} \frac{1}{cR^2} \Big(\prod_{p\leq Y} \mu_p(\Omega_p)\Big) \prod_{Y<p\leq X} \int_{\X_p} \psi(\bfx) \d\mu_p(x_p).
\end{split}
\end{equation}
Then we choose the function $\widetilde{\psi}(\bfx)$ defined by
\[
\widetilde{\psi}(\bfx) = 
\begin{cases}
1 \ \text{if} \  x_p\in\X_p & \text{for $p\leq Y$}, \\
\displaystyle\max_{s\in\K}\Big| \sum_{Y<p\leq X} \log\Big(1-\frac{x_p}{p^s} + \frac{1}{p^{2s}}\Big) \Big| \ \text{if} \ x_p\in\X_p & \text{for $Y<p\leq X$}
\end{cases}
\]
and apply \eqref{2-14} to see that
\begin{equation}
\label{3-12}
\begin{split}
\lim_{R\to\infty} &\frac{1}{cR^2} \Big(\prod_{p\leq Y} \mu_p(\Omega_p)\Big) \prod_{Y<p\leq X} \int_{\X_p} \psi(\bfx) \d\mu_p(x_p) \\
&= \Big(\prod_{p\leq Y} \mu_p(\Omega_p)\Big) \lim_{R\to\infty} \frac{1}{cR^2} \prod_{p\leq X} \int_{\X_p} \widetilde{\psi}(\bfx) \d\mu_p(x_p) \\
&=  \Big(\prod_{p\leq Y} \mu_p(\Omega_p)\Big) \lim_{R\to\infty} \frac{1}{cR^2} \sum_{r_j\leq R} \widetilde{\psi}(\bfx_j) \\
&= \Big(\prod_{p\leq Y} \mu_p(\Omega_p)\Big) \lim_{R\to\infty} \frac{1}{cR^2} \sum_{r_j\leq R} \max_{s\in\K} \Big| \sum_{Y<p\leq X} \log\Big(1-\frac{\lambda_j(p)}{p^s} + \frac{1}{p^{2s}}\Big) \Big|.
\end{split}
\end{equation}
The third assertion of Proposition~\ref{prop-2.3} follows from \eqref{3-11} and \eqref{3-12}. \fine


\bigskip
\bigskip


\bigskip
\bigskip

\noindent
Giacomo Cherubini, Dipartimento di Matematica, Universit\`a di Genova, via Dodecaneso 35, 16146 Genova, Italy. e-mail: cherubini@dima.unige.it

\noindent
Alberto Perelli, Dipartimento di Matematica, Universit\`a di Genova, via Dodecaneso 35, 16146 Genova, Italy. e-mail: perelli@dima.unige.it

\end{document}